\documentclass[11pt,english]{article}
\usepackage{amsfonts,amsmath,amsthm,amscd,amssymb,latexsym}

    \usepackage[T2A]{fontenc}
    \usepackage[cp1251]{inputenc}
    \usepackage[ukrainian]{babel}
    \usepackage{graphicx}

\begin{document}

\title{ІЗОЛЬОВАНІ СИНГУЛЯРНОСТІ ВІДОБРАЖЕНЬ З ОБЕРНЕНОЮ НЕРІВНІСТЮ ПОЛЕЦЬКОГО}

\author{Євген О.\,Севостьянов}

\theoremstyle{plain}
\newtheorem{theorem}{Теорема}[section]
\newtheorem{lemma}{Лема}[section]
\newtheorem{proposition}{Твердження}[section]
\newtheorem{corollary}{Наслідок}[section]
\theoremstyle{definition}

\newtheorem{example}{Приклад}[section]
\newtheorem{remark}{Зауваження}[section]
\newcommand{\keywords}{\textbf{Key words.  }\medskip}
\newcommand{\subjclass}{\textbf{MSC 2000. }\medskip}
\renewcommand{\abstract}{\textbf{Анотація.  }\medskip}
\numberwithin{equation}{section}

\setcounter{section}{0}
\renewcommand{\thesection}{\arabic{section}}
\newcounter{unDef}[section]
\def\theunDef{\thesection.\arabic{unDef}}
\newenvironment{definition}{\refstepcounter{unDef}\trivlist
\item[\hskip \labelsep{\bf Определение \theunDef.}]}%
{\endtrivlist}

\renewcommand{\figurename}{Мал.}

\maketitle

\begin{abstract}
Вивчаються відкриті дискретні відображення, які задовольняють вагову
оцінку спотворення модуля сімей кривих. Доведено, що ці відображення
мають неперервне продовження в ізольовану точку межі за умови
інтегровності відповідної вагової функції і належності граничної
множини відображення в даній точці до межі образу при відображенні.
\end{abstract}

\medskip {\bf Isolated singularities of mappings with the inverse Poletski
inequality.} We study open-closed discrete mappings that satisfy the
weighted estimate of the distortion of modulus of families of paths.
It is proved that the mappings mentioned above have a continuous
extension into the isolated point of the boundary, provided that the
corresponding weight function is integrable, and the cluster set of
the mapping at a given point belongs to the boundary of the image
under the mapping

\section{Вступ}
В нашій спільній роботі~\cite{SevSkv$_3$} отримано результат про
неперервне продовження в ізольовану межову точку гомеоморфізмів,
обернені до яких задовольняють вагову нерівність Полецького (див.
теорему 5.1). Основна мета даної роботи -- перенести вказаний
результат на перетворення з розгалуженням. Точніше, ми розглядаємо
відкриті дискретні відображення між двома областями розширеного
евклідового простору і припускаємо, що вони задовольняють нижню
вагову оцінку спотворення модуля сімей кривих з інтегровною
мажорантою. Зауважимо, що дослідження даної роботи знаходяться в
контексті вивчення відображень з обмеженим і скінченним спотворенням
(див., напр., \cite{GRY}--\cite{Va}). Вказані умови щодо спотворення
модуля сімей кривих добре відомі в теорії квазіконформних
відображень і їх узагальнень (див., напр.,
\cite[теорема~3.2]{MRV$_1$}, \cite[теорема~6.7.II]{Ri} і
\cite[теорема~8.5]{MRSY}).

\medskip
Звернемося до означень. Нехай $y_0\in {\Bbb R}^n,$
$0<r_1<r_2<\infty$ і
\begin{equation}\label{eq1**}
A(y_0, r_1,r_2)=\left\{ y\,\in\,{\Bbb R}^n:
r_1<|y-y_0|<r_2\right\}\,.\end{equation}
Всюди далі $\overline{{\Bbb R}^n}:={\Bbb R}^n\cup\{\infty\}.$ Для
заданих множин $E,$ $F\subset\overline{{\Bbb R}^n}$ і області
$D\subset {\Bbb R}^n$ позначимо через $\Gamma(E,F,D)$ сім'ю всіх
кривих $\gamma:[a,b]\rightarrow \overline{{\Bbb R}^n}$ таких, що
$\gamma(a)\in E,\gamma(b)\in\,F$ і $\gamma(t)\in D$ при $t \in [a,
b].$ {\it Відображенням} області $D\subset {\Bbb R}^n,$ або
$D\subset \overline{{\Bbb R}^n}$ називається довільне неперервне
перетворення $x\mapsto f(x).$  Якщо $f:D\rightarrow {\Bbb R}^n$ --
задане відображення, $y_0\in f(D)$ і
$0<r_1<r_2<d_0=\sup\limits_{y\in f(D)}|y-y_0|,$ то через
$\Gamma_f(y_0, r_1, r_2)$ ми позначимо сім'ю всіх кривих $\gamma$ в
області $D$ таких, що $f(\gamma)\in \Gamma(S(y_0, r_1), S(y_0, r_2),
A(y_0,r_1,r_2)).$ Нехай $Q:{\Bbb R}^n\rightarrow [0, \infty]$ --
вимірна за Лебегом функція.  Будемо говорити, що {\it $f$
задовольняє обернену нерівність Полецького} в точці $y_0\in f(D),$
якщо співвідношення
\begin{equation}\label{eq2*A}
M(\Gamma_f(y_0, r_1, r_2))\leqslant \int\limits_{f(D)\cap
A(y_0,r_1,r_2)} Q(y)\cdot \eta^n (|y-y_0|)\, dm(y)
\end{equation}
виконується для довільної вимірної за Лебегом функції $\eta:
(r_1,r_2)\rightarrow [0,\infty ]$ такій, що
\begin{equation}\label{eqA2}
\int\limits_{r_1}^{r_2}\eta(r)\, dr\geqslant 1\,.
\end{equation}
Зауважимо, що нерівності~(\ref{eq2*A}) добре відомі в теорії
квазірегулярних відображень і виконуються при $Q=N(f, D)\cdot K, $
де $N(f, D)$ -- максимальна кратність відображення в $D,$ а
$K\geqslant 1$ -- деяка стала, яка може бути обчислена як $K={\rm
ess \sup}\, K_O(x, f),$ $K_O(x, f)=\Vert
f^{\,\prime}(x)\Vert^n/|J(x, f)|$ при $J(x, f)\ne 0;$ $K_O(x, f)=1$
при $f^{\,\prime}(x)=0,$ і $K_O(x, f)=\infty$ при
$f^{\,\prime}(x)\ne 0,$ але $J(x, f)=0$ (див., напр.,
\cite[теорема~3.2]{MRV$_1$} або \cite[теорема~6.7.II]{Ri}).
Відображення $f:D\rightarrow {\Bbb R}^n$ називається {\it
дискретним}, якщо прообраз $\{f^{-1}\left(y\right)\}$ кожної точки
$y\,\in\,{\Bbb R}^n$ складається з ізольованих точок, і {\it
відкритим}, якщо образ будь-якої відкритої множини $U\subset D$ є
відкритою множиною в ${\Bbb R}^n.$ Як звично, покладемо
$$C(f, x):=\{y\in \overline{{\Bbb R}^n}:\exists\,x_k\in D: x_k\rightarrow x, f(x_k)
\rightarrow y, k\rightarrow\infty\}\,.$$
Тут і надалі межа $\partial D$ і замикання $\overline{D}$ області
$D$ розуміються в топології розширеного евклідового простору
$\overline{{\Bbb R}^n}.$ Виконується наступне твердження.

\medskip
\begin{theorem}\label{th2}
{Нехай $D$ і $D^{\prime}$ -- області в $\overline{{\Bbb R}^n},$
$n\geqslant 2,$ $x_0\in D,$ $f$ -- відкрите і дискретне відображення
області $D\setminus\{x_0\}$ на $D^{\,\prime},$ яке задовольняє
співвідношення~(\ref{eq2*A}) принаймні в одній скінченній точці
$y_0\in C(f, x_0).$ Нехай також $C(f, x_0)\subset
\partial D^{\,\prime}.$ Якщо $Q\in L^1(D^{\,\prime}),$ то відображення $f$ має
неперервне продовження $f\colon
D\rightarrow\overline{D^{\,\prime}}.$ Більше того, якщо
$x_0\ne\infty\ne f(x_0),$ то
для будь-якого $0<2r_0<{\rm dist}\,(x_0,
\partial D)$ і всіх $x\in
B(x_0, r_0)$ виконується нерівність
\begin{equation}\label{eq2C*}
|f(x)-f(x_0)|\leqslant\frac{C_n\cdot (\Vert
Q\Vert_1)^{1/n}}{\log^{1/n}\left(1+\frac{r_0}{|x-x_0|}\right)}
\end{equation}
де $\Vert Q\Vert_1$ -- норма функції $Q$ в $L^1(D^{\,\prime}).$}
\end{theorem}

\section{Доведення теореми~\ref{th2}}

Без обмеження загальності можна вважати $x_0\ne\infty.$ Всюди в
подальшому $h(x, y)$ позначає хордальну відстань між точками $x,
y\in \overline{{\Bbb R}^n}$ (див., напр.,
\cite[Означення~12.1]{Va}). В силу дискретності відображення $f$
існує $0<\varepsilon_0<{\rm dist}\,(x_0,
\partial D)$ таке, що $\infty\not\in f(S(x_0, \varepsilon))$ (якщо $\partial D=\varnothing,$
візьмемо довільне $\varepsilon_0>0$ зі вказаною умовою). Позначимо
$$g:=f|_{B(x_0, \varepsilon_0)\setminus\{x_0\}}\,.$$
Припустимо супротивне, а саме, що відображення $f$ не має
неперервного межового продовження в точку $x_0.$ Тоді так само і
відображення $g$ не має неперервного межового продовження в цю ж
саму точку. Оскільки простір $\overline{{\Bbb R}^n}$ є компактним,
$C(f, x_0)=C(g, x_0)\ne\varnothing.$ Тоді знайдуться $y_1, y_2\in
C(f, x_0),$ $y_1\ne y_2,$ і принаймні дві послідовності $x_m,
x^{\,\prime}_m\in B(x_0, \varepsilon_0)\setminus\{x_0\}$ такі, що
$x_m, x^{\,\prime}_m\rightarrow x_0$ при $m\rightarrow\infty,$ при
цьому, $z_m:=g(x_m)\rightarrow y_1,$
$z_m^{\,\prime}=g(x^{\,\prime}_m)\rightarrow y_2$ при
$m\rightarrow\infty.$ Можна вважати, що $y_1\ne\infty.$

\medskip
Нехай
$$D_*:=f(B(x_0, \varepsilon_0)\setminus\{x_0\})\,.$$
Покажемо, що існує $\varepsilon_1>0$ таке, що
\begin{equation}\label{eq2C}
B(y_1, \varepsilon_1)\cap f(S(x_0, \varepsilon_0))=\varnothing\,.
\end{equation}
Зауважимо, що $y_1\in \partial D_*.$ Дійсно, якщо $y_1$ -- внутрішня
точка для $D_*,$ то $y_1$ також внутрішня і для $D^{\,\prime},$
оскільки $D_*\subset D^{\,\prime}.$ Останнє суперечить умові $C(f,
x_0)\subset\partial D^{\,\prime}.$ Далі, оскільки $S(x_0,
\varepsilon_0)$ -- компакт в $D,$ то і $f(S(x_0, \varepsilon_0))$ --
компакт в $D^{\,\prime},$ тому $$h(f(S(x_0, \varepsilon_0)), C(f,
x_0))>\delta>0\,.$$ Звідси
\begin{equation}\label{eq2}
{\rm dist}\,(y_1, f(S(x_0, \varepsilon_0)))>\delta_1>0\,,
\end{equation}
де ${\rm dist}\,(A, B)$ позначає евклідову відстань між множинами
$A$ і $B$ в ${\Bbb R}^n.$ З огляду на~(\ref{eq2}),
співвідношення~(\ref{eq2C}) виконується для
$\varepsilon_1:=\delta_1.$

\begin{figure}[h]
\centerline{\includegraphics[scale=0.5]{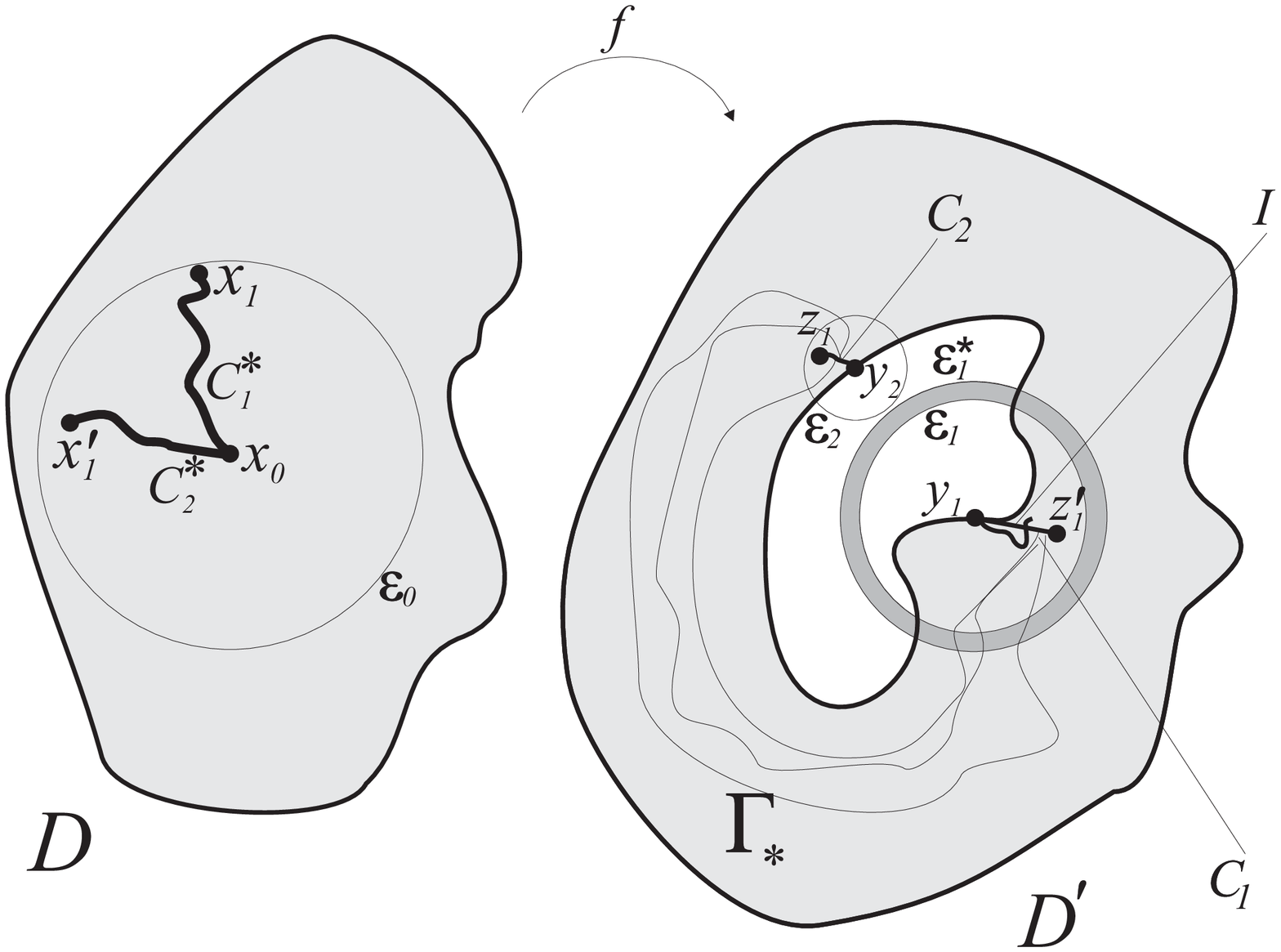}} \caption{До
доведення теореми~\ref{th2}}\label{fig1D1}
\end{figure}

\medskip
Тепер будемо міркувати наступним чином. Нехай $B_*(y_2,
\varepsilon_2)=B(y_2, \varepsilon_2)$ при $y_2\ne \infty$ і
$B_*(y_2, \varepsilon_2)=\{x\in \overline{{\Bbb R}^n}: h(x,
\infty)<\varepsilon_2\}$ при $y_2=\infty.$ Міркуючи аналогічно
доведенню співвідношення~(\ref{eq2C}), можна показати, що існує
$\varepsilon_2>0,$ таке що
\begin{equation}\label{eq5}
B_*(y_2, \varepsilon_2)\cap f(S(x_0, \varepsilon_0))=\varnothing\,.
\end{equation}
Без обмеження загальності, можна вважати, що $\overline{B(y_1,
\varepsilon_1)}\cap \overline{B_*(y_2, \varepsilon_2)}=\varnothing,$
крім того, $z_m\in B(y_1, \varepsilon_1)$ і $z^{\,\prime}_m \in
B_*(y_2, \varepsilon_2)$ при всіх $m=1,2,\ldots $ (див.
малюнок~\ref{fig1D1}).
Зауважимо, що~$B(y_1, \varepsilon_1)$ є опуклим, а, $B_*(y_2,
\varepsilon_2)$ лінійно зв'язна. В цьому випадку, точки $z_1$ і
$y_1$ можуть бути з'єднані відрізком $I(t)=z_1+t(y_1-z_1),$ $t\in
(0, 1),$ який повністю лежить в~$B(y_1, \varepsilon_1).$ Аналогічно,
точки $z^{\,\prime}_1$ і $y_2$ можна з'єднати кривою $J=J(t),$ $t\in
[0, 1],$ яка лежить в <<кулі>> $B_*(y_2, \varepsilon_2).$

\medskip
Зауважимо, що за побудовою $|I|\cap\partial D_*\ne\varnothing\ne
|J|\cap\partial D_*.$ Позначимо
$$t_*:=\sup\limits_{t\in[0, 1]: I(t)\in D_*}t\,,\qquad p_*:=\sup\limits_{t\in[0, 1]: J(t)\in D_*}t\,.$$
Нехай також
$$C_1:=I_{[0, t_*)}\,,\qquad C_2:=J_{[0, p_*)}\,.$$
За~\cite[лема~3.12]{MRV$_3$} криві $C_1$ і $C_2$ мають максимальні
підняття $C^{\,*}_1:[0, c_1)\rightarrow B(x_0,
\varepsilon_0)\setminus\{x_0\}$ і $C^{\,*}_2:[0, c_2)\rightarrow
B(x_0, \varepsilon_0)\setminus\{x_0\}$ при відображенні $g$ з
початками у точках $x_1$ і $x^{\,\prime}_1,$ відповідно. Зауважимо,
що випадок, коли $C_1(t)\rightarrow z_0$ при $t\rightarrow c_1-0,$
де $z_0\in B(x_0, \varepsilon_0)\setminus\{x_0\},$ неможливий, бо в
цій ситуації з огляду на~\cite[лема~3.12]{MRV$_3$} ми мали б, що
$c_1=t_*$ і $I(t)\rightarrow f(z_0)\in D_*,$ що суперечить означенню
$t_*.$ Тоді по~\cite[лема~3.12]{MRV$_3$}
\begin{equation}\label{eq1}
h(C^{\,*}_1(t), \partial (B(x_0,
\varepsilon_0)\setminus\{x_0\}))\rightarrow 0,\qquad t\rightarrow
c_1-0\,.
\end{equation}
Покажемо, що ситуація, коли $h(C^{\,*}_1(t), S(x_0,
\varepsilon_0))\rightarrow 0$ при $t\rightarrow c_1-0$ також є
неможливою. Дійсно, в протилежному випадку для якоїсь послідовності
$t_k\rightarrow c-0$ ми мали б, що $h(C^{\,*}_1(t_k), S(x_0,
\varepsilon_0))\rightarrow 0$ при $k\rightarrow \infty.$ В силу
компактності сфери $S(x_0, \varepsilon_0)$ знайдеться послідовність
$w_k\in S(x_0, \varepsilon_0)$ така, що $h(C_1^{\,*}(t_k), S(x_0,
\varepsilon_0))=h(C^{\,*}_1(t_k), w_k).$ Знову таки, оскільки сфера
$S(x_0, \varepsilon_0)$ компактна, то ми можемо вважати, що
$w_k\rightarrow w_0$ при $k\rightarrow \infty.$ Тоді
$C^{\,*}_1(t_k)\rightarrow w_0$ при $k\rightarrow\infty.$ Тоді за
неперервністю відображення $f$ в $D$ звідси випливає, що
\begin{equation}\label{eq7}
f(C^{\,*}_1(t_k))=C_1(t_k)\rightarrow f(w_0)\in f(S(x_0,
\varepsilon_0))
\end{equation}
при $k\rightarrow\infty.$ Останнє суперечить умові~(\ref{eq2C}), бо
одночасно $f(w_0)\in f(S(x_0, \varepsilon_0))$ і $f(w_0)\in
|I|\subset B(y_1, \varepsilon_1).$ Тоді з~(\ref{eq1}) випиває, що
\begin{equation}\label{eq1A}
h(C^{\,*}_1(t), x_0)\rightarrow 0,\qquad t\rightarrow c_1-0\,.
\end{equation}
Застосовуючи аналогічні міркування до кривої $C^{\,*}_2(t),$ можна
показати, що
\begin{equation}\label{eq1B}
h(C^{\,*}_2(t), x_0)\rightarrow 0,\qquad t\rightarrow c_2-0\,.
\end{equation}
З умов~(\ref{eq1A}) і~(\ref{eq1B}) і з огляду
на~\cite[теорема~10.12]{Va} випливає, що
\begin{equation}\label{eq1C}
M(\Gamma(|C^{\,*}_1(t)|, |C^{\,*}_2(t)|, B(x_0,
\varepsilon_0)\setminus\{x_0\}))=\infty\,.
\end{equation}
Покажемо, що~(\ref{eq1C}) суперечить умові~(\ref{eq2*A}) в точці
$y_0=y_1.$ Оскільки $\overline{B(y_1, \varepsilon_1)}\cap
\overline{B_*(y_2, \varepsilon_2)}=\varnothing,$ знайдеться
$\varepsilon^*_1>\varepsilon_1,$ для котрого ми ще маємо
$\overline{B(y_1, \varepsilon^*_1)}\cap \overline{B_*(y_2,
\varepsilon_2)}=\varnothing.$ Нехай $\Gamma_*=\Gamma(|C_1|, |C_2|,
D_*).$ Зауважимо, що
\begin{equation}\label{eq3C}
\Gamma_*>\Gamma(S(y_1, \varepsilon^*_1), S(y_1, \varepsilon_1),
A(y_1, \varepsilon_1, \varepsilon^*_1))\,.
\end{equation}
Дійсно, нехай $\gamma\in \Gamma_*,$ $\gamma:[a, b]\rightarrow {\Bbb
R}^n.$ Оскільки $\gamma(a)\in |C_1|\subset B(x_0, \varepsilon_1)$ і
$\gamma(b)\in |C_2|\subset \overline{{\Bbb R}^n}\setminus B(x_0,
\varepsilon_1),$ з огляду на~\cite[теорема~1.I.5.46]{Ku} знайдеться
$t_1\in (a, b)$ таке, що $\gamma(t_1)\in S(y_1, \varepsilon_1).$ Без
обмеження загальності, можна вважати, що
$|\gamma(t)-y_1|>\varepsilon_1$ при $t>t_1.$ Далі, оскільки
$\gamma(t_1)\in B(y_1, \varepsilon^*_1)$ і $\gamma(b)\in
|C_2|\subset {\Bbb R}^n\setminus B(x_0, \varepsilon^*_1),$ з огляду
на~\cite[Theorem~1.I.5.46]{Ku} знайдеться $t_2\in (t_1, b)$ таке, що
$\gamma(t_2)\in S(y_1, \varepsilon^*_1).$ Без обмеження загальності,
можна вважати, що $|\gamma(t)-y_1|<\varepsilon_1^*$ при $t_1<t<t_2.$
Отже, $\gamma|_{[t_1, t_2]}$ -- підкрива кривої $\gamma,$ яка
належить~$\Gamma(S(y_1, \varepsilon^*_1), S(y_1, \varepsilon_1),
A(y_1, \varepsilon_1, \varepsilon^*_1)).$ Таким чином,
співвідношення~(\ref{eq3C}) доведено.

\medskip
Встановимо тепер, що
\begin{equation}\label{eq5A}
\Gamma(|C^{\,*}_1(t)|, |C^{\,*}_2(t)|, B(x_0,
\varepsilon_0)\setminus\{x_0\})>\Gamma_f(y_1, \varepsilon_1,
\varepsilon^*_1)\,.
\end{equation}
Дійсно, якщо крива $\gamma:[a, b]\rightarrow B(x_0,
\varepsilon_0)\setminus\{x_0\}$ належить до сім'ї
$\Gamma(|C^{\,*}_1(t)|, |C^{\,*}_2(t)|, B(x_0,
\varepsilon_0)\setminus\{x_0\}),$ то $f(\gamma)$ належить $D_*,$
причому $f(\gamma(a))\in |C_1(t)|$ і $f(\gamma(a))\in |C_2(t)|,$
тобто, $f(\gamma)\in \Gamma_*.$ Тоді за доведеним вище і з огляду на
співвідношення~(\ref{eq3C}) крива $f(\gamma)$ має підкриву
$f(\gamma)^{\,*}:=f(\gamma)|_{[t_1, t_2]},$ $a\leqslant
t_1<t_2\leqslant b,$ яка належить сім'ї $\Gamma(S(y_1,
\varepsilon^*_1), S(y_1, \varepsilon_1), A(y_1, \varepsilon_1,
\varepsilon^*_1)).$ Тоді $\gamma^*:=\gamma|_{[t_1, t_2]}$ є
підкривою $\gamma$ і вона належить~$\Gamma_f(y_1, \varepsilon_1,
\varepsilon^*_1),$ що і потрібно було довести.

\medskip
Розглянемо функцію
$$\eta(t)\quad =\quad\left\{
\begin{array}{rr}
1/(\varepsilon^{*}_1-\varepsilon_1), & t\in [\varepsilon_1, \varepsilon^*_1],\\
0, & t\in {\Bbb R}\setminus[\varepsilon_1,
\varepsilon^*_1]\,.\end{array} \right.$$
Зауважимо, що $\eta$ задовольняє співвідношення~(\ref{eqA2}) при
$r_1=\varepsilon_1$ і $r_2=\varepsilon^*_1.$
Застосовуючи~(\ref{eq2*A}) в точці $y_1,$ приймаючи до уваги
умову~$Q\in L^1(D)$ і співвідношення~(\ref{eq5A}), ми отримаємо, що
$$M(\Gamma(|C^{\,*}_1(t)|, |C^{\,*}_2(t)|, B(x_0,
\varepsilon_0)\setminus\{x_0\}))\leqslant$$
\begin{equation}\label{eq4C}
\leqslant M(\Gamma_f(y_1, \varepsilon_1, \varepsilon^*_1))\leqslant
\Vert Q\Vert_1/(\varepsilon^{*}_1-\varepsilon_1)^{n} <\infty\,,
\end{equation}
де $\Vert Q\Vert_1$ позначає $L^1$-норму функції $Q$ в області
$D^{\,\prime}.$
Співвідношення~(\ref{eq1C}) і~(\ref{eq4C}) суперечать одне одному,
що і спростовує припущення про наявність різних $y_1$ і $y_2\in C(f,
x_0).$

\medskip
Нарешті, якщо $x_0\ne \infty,$ то розглянемо область
$D_1:=D\setminus\{f^{\,-1}(\infty)\}.$ Зауважимо, що в силу
замкненості відображення $f$ множина $\{f^{\,-1}(\infty)\}$ є
скінченною, див.~\cite[лема~3.3]{MS}, тому $D_1$ є областю, а точка
$x_0$ є її внутрішньою точкою. В такому випадку,
співвідношення~(\ref{eq2C*}) є результатом
роботи~\cite[Теорема~1.1]{SSD}.~$\Box$

КОНТАКТНА ІНФОРМАЦІЯ

\medskip
\noindent{{\bf Євген Олександрович Севостьянов} \\
{\bf 1.} Житомирський державний університет ім.\ І.~Франко\\
кафедра математичного аналізу, вул. Велика Бердичівська, 40 \\
м.~Житомир, Україна, 10 008 \\
{\bf 2.} Інститут прикладної математики і механіки
НАН України, \\
вул.~Добровольського, 1 \\
м.~Слов'янськ, Україна, 84 100\\
e-mail: esevostyanov2009@gmail.com}


\begin{thebibliography}{99}

\bibitem{SevSkv$_3$} Sevost’yanov, E.A. and Skvortsov~S.A.,
\emph{On mappings whose inverse satisfy the Poletsky inequality},
Ann. Acad. Scie. Fenn. Math., \textbf{45}, 2020, p.~259--277.

\bibitem{GRY} Gutlyanskii~V.~Ya., Ryazanov~V.~I.,
Yakubov~E., \emph{The Beltrami equations and prime ends},
Український математичний вiсник, \textbf{12}, № 1, 2015, p. 27–-66;
translation \emph{The Beltrami equations and prime ends}, J. Math.
Sci. (N.Y.), \textbf{210}, no.~1, 2015, p.~22–-51.

\bibitem{MRV$_1$} Martio~O., Rickman~S., and V\"{a}is\"{a}l\"{a}~J.,
\emph{Definitions for quasiregular mappings}, Ann. Acad. Sci. Fenn.
Ser. A1., \textbf{448}, 1969, p.~1--40.

\bibitem{MRV$_3$} Martio~O., Rickman~S., and V\"{a}is\"{a}l\"{a}~J.,
\emph{Topological and metric properties of quasiregular mappings},
Ann. Acad. Sci. Fenn. Ser. A1., \textbf{488,} 1971, p.~1--31.

\bibitem{MRSY} Martio O., Ryazanov V., Srebro U. and Yakubov
E., \emph{Moduli in Modern Mapping Theory.} -- New York: Springer
Science + Business Media, LLC, 2009.

\bibitem{MS} Martio~O., Srebro~U., \emph{Periodic quasimeromorphic
mappings in ${\Bbb R}^n$}, J. d'Anal. Math. \textbf{28}, no.~1,
1975, p.~20--40.

\bibitem{Ri} Rickman S., \emph{Quasiregular mappings.} -- Berlin: Springer-Verlag, 1993.

\bibitem{Va} V\"{a}is\"{a}l\"{a} J., \emph{Lectures on $n$-Dimensional
Quasiconformal Mappings}, Lecture Notes in Math. \textbf{229},
Berlin etc.: Springer--Verlag, 1971.

\bibitem{Ku} Куратовский~К., \emph{Топология}, т.~2. -- М.:
Мир, 1969.

\bibitem{SSD} Sevost'yanov~E.A., Skvortsov~S.O., Dovhopiatyi~O.P.,
\emph{On mappings satisfying the inverse Poletsky inequality} //
www. arxiv. org, arXiv:1904.01513.

\end{thebibliography}
\end{document}